# A timestepper approach for the systematic bifurcation and stability analysis of polymer extrusion dynamics

M.E. Kavousanakis[1], L. Russo[2], C. I. Siettos[3,*], A. G. Boudouvis[1], G.C. Georgiou[4]

[1]School of Chemical Engineering, National Technical University of Athens,
GR 157 80, Athens, Greece

[2]Department of Chemical Engineering, University "Federico II" of Naples,
80125, Piazzale Tecchio 80, Naples, Italy

[3]School of Applied Mathematics and Physical Sciences,
National Technical University of Athens, GR 157 80, Athens, Greece

[4]Department of Mathematics and Statistics, University of Cyprus, PO Box 20537, 1678,
Nicosia, Cyprus

**Abstract**

We discuss how matrix-free/timestepper algorithms can efficiently be used with dynamic non-Newtonian fluid mechanics simulators in performing systematic stability/bifurcation analysis. The timestepper approach to bifurcation analysis of large scale systems is applied to the plane Poiseuille flow of an Oldroyd-B fluid with non-monotonic slip at the wall, in order to further investigate a mechanism of extrusion instability based on the combination of viscoelasticity and nonmonotonic slip. Due to the nonmonotonicity of the slip equation the resulting steady-state flow curve is nonmonotonic and unstable steady-states appear in the negative-slope regime. It has been known that self-sustained oscillations of the pressure gradient are obtained when an unstable steady-state is perturbed [Fyrillas et al., Polymer Eng. Sci. 39 (1999) 2498-2504].

Treating the simulator of a distributed parameter model describing the dynamics of the above flow as an input-output black-box timestepper of the state variables, stable and unstable branches of both equilibrium and periodic oscillating solutions are computed and their stability is examined. It is shown for the first time how equilibrium solutions lose stability to oscillating ones through a subcritical Hopf bifurcation point which generates a branch of unstable limit cycles and how the stable periodic solutions lose their stability through a critical point which marks the onset of the unstable limit cycles. This implicates the coexistence of stable equilibria with stable and unstable periodic solutions in a narrow range of volumetric flow rates.

---

[*] Corresponding author

# 1. Introduction

The complex viscoelastic character of polymers, the normal stress differences and the high extensional viscosity and wall slip may lead to nonlinear phenomena and undesirable instabilities in polymer processing. Time-periodic phenomena are often observed, such as pressure oscillations at fixed volumetric flow rate in the stick-slip extrusion instability and draw resonance, which gives rise to spontaneous thickness and width oscillations in film casting, to a periodic fluctuation of the cross-sectional area in fiber spinning, and to periodic fluctuations of the bubble diameter in film blowing [1]. For example, draw resonance in the latter process corresponds to a self-sustained limit cycle type supercritical Hopf bifurcation [2]. In such flows, in addition to the steady-state solutions and linear stability, transient studies and nonlinear stability analyses are necessary, in order to develop techniques for process optimization [2].

Modelling and understanding the mechanisms of such flow instabilities by determining the regions and the critical points where these occur are of major importance. For this purpose, the efficient simulation of the transient behaviour of the underlying physical system, usually expressed in terms of a system of ordinary differential and algebraic equations or integro-partial differential algebraic equations, is required. Over the last years, some excellent temporal (direct integration in time) commercial and home-made fluid mechanics simulation packages are the tools of choice. Such packages, which incorporate many man-years of effort and expertise, may also allow the computation of steady states using Newton-like solvers and comprise a very good option even for large scale systems.

However, other important tasks, such as the exact location of the critical points that mark the onset of instabilities, as well as the dependence of the location on parameter values, cannot be obtained easily using temporal simulations. Furthermore, the tracing of branches of unstable steady states is often (in the absence of good initial guesses) impossible for large-scale systems. For the systematic and accurate analysis of the model dynamics, one has to resort to bifurcation analysis. Numerical bifurcation theory provides an arsenal of algorithms and software packages, such as AUTO, CONTENT and MATCONT for tasks such as the continuation of stable or unstable steady states and limit cycles, and the continuation of critical points [3-6]. While, these software packages are invaluable tools for performing systematic analysis for small to medium scale systems, there are some drawbacks in using them. Most of them require as input the system evolution equations, which are assumed to be explicitly available in discretized form. Linking the evolution equations with such packages is not a trivial task. These packages often use a Newton-like method, which requires the calculation of the Jacobian of the system (i.e. the matrix with the partial derivatives of the "right-hand-sides" of the discretized governing ODEs or the PDEs, with respect to the discrete unknowns). This imposes a serious computational burden in the analysis of large-scale systems. But even for small- to medium-size systems, tasks, such as the continuation of limit cycles or the continuation of turning points of limit cycles, become overwhelmingly computationally expensive or even prohibitive; these computations are usually performed by augmenting the system space with one more variable corresponding to the normalized time variable. The latter turns an initial-value problem to a boundary-value one, with a consequent vast increase in the size of the problem.

The solution to the "curse of dimensionality" comes from the matrix-free methods of iterative linear algebra [7]. Here one does not need to numerically compute the Jacobian of the system. What is required is the calculation of matrix-vector products for a sequence of known vectors in order to perform tasks such as solving (with the Newton method) large systems of non-linear equations or finding the critical-leading eigenvalues. Typically the dynamic simulator is only treated as a black box input-output timestepper which resolves the linking problem. The necessary information is acquired by Newton-Picard schemes, such as the Recursive Projection Method (RPM) of Shroff and Keller [8] "wrapped around" the commercial simulator. "Wrapping" of matrix-free algorithms around industrial process simulators (like gPROMS) has been recently discussed by Siettos et al. [9] who applied the RPM for the efficient location of the cycle steady states and stability analysis of a periodically forced process.

The purpose of this paper is twofold: (a) to introduce the concept of matrix-free/timestepper approach that enables non-Newtonian fluid dynamics simulators to perform efficiently numerical stability/bifurcation analysis (such as continuation of both steady states and limit cycles and the computation of their stability); and (b) to demonstrate the applicability of the method to viscoelastic flow problems in performing systematic numerical bifurcation and stability analysis of periodic solutions.

For the latter objective, the time-dependent, one-dimensional plane Poiseuille flow of an Oldroyd-B fluid with nonmonotonic slip at the wall has been chosen. This problem has been considered by Georgiou and co-workers [10-12] who studied the combined effect of elasticity and non-monotonic slip and examined whether this can provide an explanation for the stick-slip extrusion instability [1]. All theoretical explanations suggested in the literature for this instability are based on the non-monotonicity of the flow curve (the plot of the wall shear stress versus the apparent shear rate, or, equivalently, the plot of the pressure gradient versus the volumetric flow rate), which exhibits a maximum and a minimum, and the fact that steady-state solutions corresponding to the negative slope regime of the flow curve are unstable [13] The transitions from the maximum of the flow curve to the right positive-slope branch and from the minimum to the left positive slope branch lead to a limit cycle, which describes the observed pressure and flow-rate oscillations in flow-rate-controlled experiments [14].

Non-monotonicity of the flow curve can be obtained by a non-monotonic slip law (adhesive failure) or by a non-monotone constitutive equation (bulk failure). In the proposed explanations involving slip, this is combined with either compressibility or elasticity. The important role of slip in the stick-slip instability, indicated by both indirect and direct wall slip measurements, is also supported by the fact that only mechanisms involving slip lead to self-sustained pressure oscillations and generate waves on the extrudate surface. However, in addition to the experimental evidence for the importance of the compressibility of the melt in the reservoir, only the compressibility/slip mechanism can lead to persistent pressure and flow rate oscillations between the two stable branches of the flow curve. The periodic transitions between a weak slip (or no-slip) and a strong slip at the capillary wall (i.e. the jumps between the two branches of the flow curve) which lead to the pressure and flow rate oscillations are sustained by the compressibility of the melt in the reservoir. This mechanism has been employed in various one-dimensional phenomenological models describing the stick-slip instability (see [15] and references therein) and in two-dimensional simulations of both Poiseuille and extrudate-swell flows (see [14] and references therein).

Viscoelasticity may replace compressibility and, when combined with non-monotonic slip, can act as a storage of elastic energy generating self-sustained pressure oscillations and waves on the extrudate surface in the stick-slip regime. Due to the absence of compressibility, however, this mechanism cannot generate jumps of the volumetric flow rate between the two stable branches of the flow curve and leads only to small-amplitude small-wavelength distortions of the extrudate surface consistent with sharkskin rather than with the stick-slip instability. Nevertheless, these may be superimposed to the much larger oscillations caused by the compressibility/slip mechanism.

Georgiou and co-workers [10-12] explored the elasticity/slip mechanism by solving numerically the axisymmetric Poiseuille and extrudate-swell flows of an Oldroyd-B fluid with non-monotonic slip at the wall. The stability of the one-dimensional Poiseuille flow was investigated by means of linear stability analysis and numerical time-dependent calculations, which revealed the existence of unstable steady-state solutions in the negative-slope regime of the flow curve. In the latter regime, self-sustained pressure oscillations are observed, the amplitude and the period of which increase with elasticity [10-12]. The numerical solutions of the two-dimensional axisymmetric Poiseuille and extrudate-swell flows showed the existence of a second mode of periodicity in the axial direction and the generation of small-amplitude waves on the extrudate surface [11].

Shore and co-workers [16-18] presented a phenomenological hydrodynamic model describing the Poiseuille flow of a Maxwell fluid with slip at the wall. They assumed that the conformation of polymers near the surface undergoes a first-order transition as a function of the wall shear stress. This conformational change produces stick-slip behavior and leads to a multivalued shear stress/slip velocity curve. Using linear stability analysis as well as numerical calculations of the linearized, incompressible, two-dimensional momentum equations, they demonstrated the existence of self-sustained oscillations, and related the period of the latter to the sharkskin polymer extrusion instability.

Black and Graham [19, 20] demonstrated that the combination of elasticity with a *monotonic* slip model that takes into account the normal stress dependence of the slip velocity leads to short wavelength shear flow instabilities at sufficiently high shear rates, which is also qualitatively consistent with experimental observations of the sharkskin instability. They considered the plane simple shear and Poiseuille flows of both the upper-convected Maxwell and Phan-Thien-Tanner fluids and reported that the scaling of the critical shear stress for instability with modulus and molecular weight and of the distortion period with polymer relaxation time are qualitatively consistent with experimental observations of the sharkskin instability in linear polyethylenes.

The paper is organized as follows. We start by presenting the timestepper approach for the continuation of solution branches of limit cycles for large scale systems. In Section 3, the equations governing the plane Poiseuille flow of an Oldroyd-B fluid along with the nonmonotonic slip equation are provided and the steady-state solutions are discussed. In Section 4 the results of the bifurcation analysis obtained through the timestepper approach are presented and discussed. It is shown for the first time that the loss of stability of steady-state solutions to sustained oscillating ones takes place through a subcritical Hopf bifurcation while the branch of stable periodic solutions loses stability at a critical point of limit cycles. This combination implicates the coexistence of stable steady states with stable and unstable periodic solutions in a narrow range of values of flow rates that can drive the system to abrupt loss of stability. Finally, Section 5 summarizes the conclusions.

## 2. Computations of periodic solutions with the timestepper approach

Consider the continuous time, autonomous nonlinear system

$$\frac{d\boldsymbol{u}}{dt} = \boldsymbol{F}(\boldsymbol{u},\mu), \tag{1}$$

where $\boldsymbol{u}$ denotes the state vector, accessible through measurement, $\mu$ is a system parameter serving here as the bifurcation parameter, and $\boldsymbol{F}$ is a not explicitly available function. The question is how to systematically construct bifurcation diagrams of steady and periodic solutions and perform their stability analysis. The answer comes from the concept of timestepping [8, 9, 21-24]. Despite the fact that $\boldsymbol{F}$ is not explicitly available, we still assume that we do have a discrete input-output map (a "black box" timestepper) that, given the initial state of the system $(\boldsymbol{u}(0), \mu)$ reports the solution, the result of the integration, after a given time horizon $T_h$, i.e.,

$$\boldsymbol{u}(T_h) = \boldsymbol{\Phi}_{T_h}(\boldsymbol{u}(0), \mu), \tag{2}$$

where $\boldsymbol{\Phi}_{T_h}$ is the temporal evolution operator of the system. At an outer level, and depending on the task to be carried out (such as the location of fixed points, the design of a feedback controller and optimisation), well established numerical analysis algorithms (such as the Newton-Raphson technique) can be utilized to *estimate* "on demand" the required quantities, e.g. residuals, Jacobians and Hessians [25, 26]. The idea is extremely simple: the timestepper is called as a black-box subroutine with nearby appropriately perturbed initial conditions and for relatively short time intervals.

If a periodic oscillatory behavior is observed then one seeks for solutions which satisfy

$$\boldsymbol{u}(0) = \boldsymbol{u}(T), \tag{3}$$

with $T$ denoting the period of oscillation. Using $T$ as sampling time, periodic solutions can be computed as fixed points of the mapping

$$\boldsymbol{\Phi}_T : \boldsymbol{u} - \boldsymbol{\Phi}_T(\boldsymbol{u}, \mu) = \boldsymbol{0} \tag{4}$$

augmented by the so-called phase constraint (also called a pinning condition)

$$g(\boldsymbol{u}, \mu, T) = 0, \tag{5}$$

which factors out the infinite members of the family of periodic solutions. The above condition enables the computation of the unknown period $T$.

The tracing of the branches of periodic solutions can be achieved using the linearized pseudo arc-length condition [3]:

$$N(\boldsymbol{u}, \mu, T) = \boldsymbol{\alpha} \cdot (\boldsymbol{u} - \boldsymbol{u}_1) + \beta(\mu - \mu_1) + \gamma(T - T_1) - \delta s = 0, \tag{6a}$$

where

$$\boldsymbol{\alpha} \equiv \frac{(\boldsymbol{u}_1 - \boldsymbol{u}_0)^T}{\delta s}, \quad \beta \equiv \frac{(\mu_1 - \mu_0)}{\delta s}, \quad \gamma \equiv \frac{(T_1 - T_0)}{\delta s} \tag{6b}$$

and $\delta s$ is the pseudo arc-length continuation step. $(\boldsymbol{u}_0, \mu_0, T_0)$ and $(\boldsymbol{u}_1, \mu_1, T_1)$, represent two already known periodic solutions. Equation (6) constrains the "next" periodic solution to lie on a hyperplane perpendicular to the tangent of the bifurcation diagram at $(\boldsymbol{u}_1, \mu_1, T_1)$, approximated through $(\boldsymbol{\alpha}, \beta, \gamma)$, at a distance δs from it. The computation of the periodic solutions can be obtained using an iterative procedure like the Newton-Raphson technique [7]. The procedure involves the iterative solution of the following linearized system:

$$\begin{bmatrix} I - \dfrac{\partial \boldsymbol{\Phi}_T}{\partial \boldsymbol{u}} & -\dfrac{\partial \boldsymbol{\Phi}_T}{\partial \mu} & -\dfrac{\partial \boldsymbol{\Phi}_T}{\partial T} \\ \dfrac{\partial g}{\partial \boldsymbol{u}} & \dfrac{\partial g}{\partial \mu} & \dfrac{\partial g}{\partial T} \\ \boldsymbol{a} & \beta & \gamma \end{bmatrix} \begin{bmatrix} \delta \boldsymbol{u} \\ \delta \mu \\ \delta T \end{bmatrix} = - \begin{bmatrix} \boldsymbol{u} - \boldsymbol{\Phi}_T(\boldsymbol{u}, \mu) \\ g(\boldsymbol{u}, \mu, T) \\ N(\boldsymbol{u}, \mu, T) \end{bmatrix} \qquad (7)$$

Note that for the computation of the sub-Jacobian $\dfrac{\partial \boldsymbol{\Phi}_T}{\partial \boldsymbol{u}}$ and the derivatives $\dfrac{\partial \boldsymbol{\Phi}_T}{\partial \mu}$, $\dfrac{\partial \boldsymbol{\Phi}_T}{\partial T}$, no explicit evolution equations are required. Their approximation can be achieved numerically by calling the timestepper at appropriately perturbed values of the corresponding unknowns. The above framework enables the temporal simulator to converge to both stable and unstable periodic solutions and trace their locations, i.e. to fulfill tasks that the simulator was not explicitly designed for [23].

However, the above approach turns to be computationally inefficient especially for large scale problems. An alternative approach is to solve (4) and (5) at each Newton step with a matrix-free iterative solver, such as the Newton-Generalized Minimum Residual (Newton-GMRES) method [7] (see also Fig. 1). The advantage of using GMRES is that the explicit calculation and storage of the Jacobian is not required. Only matrix-vector multiplications are needed which can be performed at low cost by calling the timestepper from *nearby* initial conditions allowing the estimation of the action of the linearization of a map $\boldsymbol{\Phi}_{T_h}$ on known vectors, since

$$D\boldsymbol{\Phi}_{T_h}(\boldsymbol{u}) \cdot \boldsymbol{q} \approx \dfrac{\boldsymbol{\Phi}_{T_h}(\boldsymbol{u} + \varepsilon \boldsymbol{q}) - \boldsymbol{\Phi}_{T_h}(\boldsymbol{u})}{\varepsilon}, \qquad (8)$$

where $\varepsilon$ is a small and appropriately chosen scalar. At step $j$ the algorithm minimizes the residual

$$\boldsymbol{R} \equiv \boldsymbol{u} - \boldsymbol{\Phi}_{T_h}(\boldsymbol{u}, \mu) \qquad (9)$$

by producing an orthonormal basis $\{q_1, q_2, ...., q_j\}$ of the Krylov subspace $K_j$ spanned by $\{q_1, D\boldsymbol{\Phi}_{T_h}(\boldsymbol{u}) \cdot q_1, ...., D\boldsymbol{\Phi}_{T_h}(\boldsymbol{u})^{j-1} \cdot q_{j-1}\}$. The projection of $D\boldsymbol{\Phi}_{T_h}(\boldsymbol{u})$ on $K_j$ is represented in the basis $V_j \equiv \{q_j\}$ by the upper Hessenberg matrix

$$H_j = V_j^T D\boldsymbol{\Phi}_{T_h}(\boldsymbol{u}) V_j \qquad (10)$$

whose elements are the coefficients $h_{ij}$.

Alternative algorithms, such as the Recursive Projection Method of Shroff and Keller [8], and other Newton-Picard methods, such as the ones presented by Lust et al. [21], can also be used to compute in an efficient manner both steady states and periodic solutions and construct their bifurcation diagrams. These algorithms use the timestepper to approximate iteratively the slow stable or unstable eigenspace of $D\boldsymbol{\Phi}_{T_h}(\boldsymbol{u})$ which is assumed to be of low dimension. Newton's method is implemented on the slow eigenspace to accelerate the convergence to the fixed point, while Picard iteration (time integration) is applied on the complement of that subspace to eliminate the fast dynamics of the mapping given by Eq. (2).

Having calculated the periodic solution and its period *T*, the dynamic behavior of the particular system can be described by the stroboscopic/Poincaré map: the system trajectories are not recorded continuously in time, but only once every period. Under this representation solutions appear stroboscopically as fixed points of the mapping $\boldsymbol{\Phi}_T$ under which $\boldsymbol{u}^* = \boldsymbol{\Phi}_T(\boldsymbol{u}^*, \mu)$. For a fixed value of the bifurcation parameter $\mu$ the equation of motion in a neighborhood of a fixed point can be described by the linear Poincaré mapping: $\tilde{\boldsymbol{u}}_{k+1} = D\boldsymbol{\Phi}_T(\boldsymbol{u})\tilde{\boldsymbol{u}}_k$, where $\tilde{\boldsymbol{u}} = \boldsymbol{u} - \boldsymbol{u}^*$. The eigenvalues (Floquet multipliers) of the matrix $D\boldsymbol{\Phi}_T$ (the monodromy matrix) determine the local stability of the system. Again, the leading Floquet multipliers can be estimated using an ad hoc matrix-free iterative eigensolver such as the Arnoldi procedure [27, 28] or delivered as a byproduct of the Newton-Picard iterations.

## 3. Plane Poiseuille flow of an Oldroyd-B fluid

In the Oldroyd-B model the extra stress tensor **T** is decomposed into a viscoelastic part, **T₁**, and a purely viscous part, as follows:

$$\mathbf{T} = \mathbf{T_1} + \eta_2 \left[ \nabla \mathbf{v} + (\nabla \mathbf{v})^T \right] \tag{11}$$

where $\eta_2$ is a viscosity parameter, **v** is the velocity vector, $\nabla \mathbf{v}$ is the velocity-gradient tensor, and (*T*) denotes the transpose (not to be confused with the period in the previous section). The tensor **T₁** is defined by

$$\mathbf{T_1} + \lambda \left[ \frac{\partial \mathbf{T_1}}{\partial t} + \mathbf{v} \cdot \nabla \mathbf{T_1} - (\nabla \mathbf{v})^T \cdot \mathbf{T_1} - \mathbf{T_1} \cdot \nabla \mathbf{v} \right] = \eta_1 \left[ \nabla \mathbf{v} + (\nabla \mathbf{v})^T \right], \tag{12}$$

where $\lambda$ is the relaxation time and $\eta_1$ is another viscosity parameter. The shear viscosity is given by $\eta_1 + \eta_2$

We consider the one-dimensional plane Poiseuille flow of an Oldroyd-B fluid with slip at the wall. The geometry of the flow is shown in Fig. 2. Scaling the lengths by the channel half-width H, the velocity vector by a characteristic velocity V, the pressure and the stress components by ($\eta_1 + \eta_2$)V/H, and the time by H/V, the dimensionless x-momentum equation takes the form

$$Re \frac{\partial v_x}{\partial t} = -\nabla P + \frac{\partial T_1^{xy}}{\partial y} + \eta_2 \frac{\partial^2 v_x}{\partial y^2}, \tag{13}$$

where $-\nabla P$ denotes the dimensionless pressure gradient, *Re* is the Reynolds number,

$$Re = \frac{\rho V H}{\eta_1 + \eta_2} \tag{14}$$

and $\rho$ is the density. The only nonzero component of **T₁** is the xy-component for which one gets:

$$T_1^{xy} + We \frac{\partial T_1^{xy}}{\partial t} = (1 - \eta_2) \frac{\partial v_x}{\partial y}. \tag{15}$$

where We is the Weissenberg number defined by

$$We = \frac{\lambda V}{H} \tag{16}$$

As shown in Fig. 2, along the symmetry plane ($y=0$) the velocity gradient is zero. Along the wall ($y=1$), slip is assumed to occur following a slip law of the following general form

$$\sigma_w = F(v_w) \tag{17}$$

where $\sigma_w$ is the wall shear stress,

$$\sigma_w = \left( T_1^{xy} + \eta_2 \frac{\partial v_x}{\partial y} \right)\bigg|_{y=1}$$

and $v_w$ is the slip velocity.

In steady state, the solution of the system of Eqs. (13) and (15) is:

$$v_x(y) = v_w - \frac{1}{2}\nabla P(1 - y^2) \tag{18}$$

and

$$T_1^{xy}(y) = \eta_1 \nabla P\, y \tag{19}$$

so that

$$T^{xy} = T_1^{xy} + \eta_2 \frac{dv_x}{dy} = \nabla P\, y \tag{20}$$

(the viscosities $\eta_1$ and $\eta_2$ are dimensionless and $\eta_1 + \eta_2 = 1$). The slip velocity $v_w$ satisfies the condition

$$\sigma_w = -\nabla P = F(v_w) \tag{21}$$

and the volumetric flow rate Q is given by

$$Q = v_w - \frac{1}{3}\nabla P = v_w + \frac{1}{3}F(v_w) \tag{22}$$

Let us consider here the arbitrary non-monotonic slip equation used by Georgiou and co-workers [10, 29] and others [15]):

$$\sigma_w = F(v_w) = A_1\left(1 + \frac{A_2}{1 + A_3 v_w^2}\right) v_w . \tag{23}$$

where $A_1$, $A_2$ and $A_3$ are dimensionless parameters. The nonmonotonic flow curve for $A_1=1$, $A_2=15$ and $A_3=100$ is shown in Fig. 3.

As noted by Fyrillas et al. [12], in the case of time-dependent flow with fixed volumetric flow rate, the integration of Eq. (13) over the channel cross section reveals that

$$-\nabla P(t) = F(v_w(t)) \quad \forall t \tag{24}$$

In other words, when it is plotted as a function of $v_w$, the time-dependent pressure gradient follows the graph of $F(v_w)$. The linear stability of the steady-state solutions to infinitesimal one-dimensional antisymmetric disturbances at fixed volumetric flow rate has been studied by Fyrillas et al. [12]. The stability of the basic solutions has been found to depend on three parameters: $\eta_2$, $F'(v_w^*)$, where $v_w^*$ is the steady-state slip velocity, and the elasticity number $E = We/Re$. Unstable solutions exist only when $F'(v_w^*)$ is negative and above the marginal stability curve of $-F'(v_w)$ versus $E$ for a given value of $\eta_2$. The flow is stable for a Newtonian fluid ($\eta_2 = 1$) and destabilizes as $\eta_2$ is decreased and/or $E$ is increased [12].

## 4. Timestepper based bifurcation and stability analysis

Second-order central finite differences are used for the numerical approximation of the $v_x$ and $T_1^{xy}$ spatial derivatives (partial derivatives with respect to $y$) and a fully-implicit (Euler backward difference) scheme is applied for time integration. The one dimensional domain is discretized with $N=801$ equidistant nodes and the time step is fixed to $dt=10^{-5}$. In the computations we use the same set of parameters as that employed by Fyrillas et al. [12]: $Re=1$, $We=0.1$, $\eta_2=0.1$, $A_1=1$, $A_2=15$ and $A_3=100$.

The determination of critical points which initiate regions of instability was performed with timestepper based linear stability analysis. In order to characterize the stability of a steady-state solution $\boldsymbol{u}^*$, the critical eigenvalues of the steady-state Jacobian matrix, i.e. the eigenvalues with maximum real part, are sought. The steady state solution satisfies the relation:

$$\boldsymbol{R}(\boldsymbol{u}^*, \mu) \equiv \boldsymbol{u}^* - \boldsymbol{\Phi}_{T_h}(\boldsymbol{u}^*, \mu) = 0 \tag{25}$$

where $T_h$ is the time horizon of the "black-box" timestepper. The eigenvalues $\sigma_i$ of the matrix

$$\left.\frac{\partial \boldsymbol{R}}{\partial \boldsymbol{u}}\right|_{(\boldsymbol{u}^*, \mu)} = \boldsymbol{I} - \left.\frac{\partial \boldsymbol{\Phi}_{T_h}}{\partial \boldsymbol{u}}\right|_{(\boldsymbol{u}^*, \mu)} \tag{26}$$

where $\boldsymbol{I}$ is the identity matrix, can be obtained from the application of the Arnoldi eigensolver using MATLAB's *eigs* function. Arnoldi's algorithm operates in a matrix-free context which can be used to estimate the most critical eigenvalues. The eigenvalues $\kappa_i$ of the Jacobian matrix $\left.\frac{\partial \boldsymbol{\Phi}_{T_h}}{\partial \boldsymbol{u}}\right|_{(\boldsymbol{u}^*, \mu)}$ of the discrete map (2) are related to the eigenvalues $\lambda_i$ of the Jacobian matrix $\left.\frac{\partial \boldsymbol{F}}{\partial \boldsymbol{u}}\right|_{(\boldsymbol{u}^*, \mu)}$ of the continuous problem (1) as follows [23]:

$$\kappa_i = exp(\lambda_i T_h), \tag{27}$$

A steady-state solution is stable when all eigenvalues $\lambda_i$ have negative real part and unstable when at least one eigenvalue has positive real part. In our computation the time horizon of the time stepper was chosen to be $T_h$=10$^{-3}$ while in our difference approximation of the matrix-product operations we used perturbations of the order of 10$^{-6}$.

The application of linear stability analysis to the steady-state solutions indicated the existence of a Hopf point at a volumetric flow rate Q≈0.4135 from which a branch of stable limit cycle solutions emerges. The critical eigenvalues $\lambda_i$ of the Jacobian matrix are depicted in Fig. 4, where the transition of two complex eigenvalues from the left real plane to the right plane can be seen. A second, Hopf point was computed at $Q$≈0.5215 (see Fig. 5), from which a branch of unstable limit cycles emerges, as will be discussed below. It should be pointed out that the two Hopf points correspond to the same critical value of $-F'(v_w^*)$ (i.e., about 0.384) at which the flow for $\eta_2$=1 and $E$=10 becomes linearly unstable. The timestepper approach calculations agree well with the linear stability analysis results provided by Fyrillas et al. [12]. The basic difference is that in the present work the stability analysis is carried out numerically by exploiting the numerical integrator of the governing system of equations.

The existence of periodic solutions at the parametric region specified by the two Hopf points can also be proven by executing transient simulations of slightly perturbed (unstable) steady state solutions [12]. In Fig. 6, we show the evolution of the slip velocity and the pressure gradient when the steady-state solution computed at Q=0.449 is perturbed by setting Q=0.45. The solution becomes practically periodic after one oscillation. It is clear that elasticity acts as the storage of elastic energy leading to self-sustained oscillations in the unstable regime (i.e. the region between the two Hopf points). The abrupt changes in the slip velocity correspond to jumps from weak to strong slip and vice versa. As discussed in Section 3, the time-dependent pressure gradient oscillations along the graph of the slip function $F(v_w^*)$. The local minima are due to the fact that as $-\nabla P(t)$ oscillates to the left, it passes not only the left Hopf point but also the maximum of the slip function and it thus decreases slightly reaching a local minimum after which it starts moving to the right. The explanation for the local maxima is similar. Other transient simulations have also been carried out in the neighborhoods of the Hopf points where the perturbed unstable steady-state solutions demonstrate initially an oscillatory behavior which ultimately converges to (stable) limit cycle solutions.

In order to perform a more comprehensive study of the periodic solutions which "live" in the parametric region marked by the two Hopf points, the periodic solution branch is traced by applying a pseudo arc-length continuation algorithm, as described in Section 2. The continuation algorithm enables the possibility of tracing both stable and unstable periodic solutions. The volumetric flow rate Q is treated as the continuation parameter and in our computations the phase constraint (5) has the following form [21]:

$$\left.\frac{\partial \boldsymbol{u}^{ref}}{\partial t}\right|_{t=0} \cdot \left(\boldsymbol{u}-\boldsymbol{u}^{ref}\right)=0. \tag{28}$$

This condition minimizes the phase shift of the sought solution $\boldsymbol{u}$ with respect to a reference solution $\boldsymbol{u}^{ref}$, e.g. the computed periodic solution at a previous parameter value.

The complete bifurcation diagram of steady-state and periodic solutions (stable and unstable) is depicted in Fig. 7, where the pressure gradient is plotted as a function of the volumetric flow rate Q. The solid line corresponds to stable steady states, the dashed one to unstable steady states, the solid line with full circles to stable periodic solutions, and the dashed line with open circles to unstable periodic solutions. The circles correspond to the maximum value of the periodic pressure gradient. A branch of stable periodic solutions is initiated at the first Hopf point (Q≈0.4135). The parametric study revealed that stable periodic solutions also exist beyond the right Hopf point (Q≈0.5125) up to a critical turning point at Q≈0.5364 which marks the onset of unstable periodic solutions. The local stability of the periodic solutions was determined from the Floquet multipliers of the monodromy matrix $D\boldsymbol{\Phi}_T$, which were computed by a matrix-free eigensolver (Arnoldi). As predicted from the theory, there is always a Floquet multiplier equal to one, which corresponds to a phase-shift along the periodic solution. The periodic solution is stable if all the other Floquet multipliers have modulus smaller than one and unstable if there exists at least one Floquet multiplier with modulus greater than one. In Fig. 8 the dominant Floquet multipliers are shown and one can see the transition of the leading eigenvalue in the outer region specified by the complex unit circle.

The computed periodic solutions lose their stability at Q≈0.5364. A focus on the turning point of periodic solutions is shown at the inset of Fig. 7. The emerging (relatively) steep branch of unstable periodic solutions joins the steady-state solution branch at the second subcritical Hopf point (Q≈0.5215). This implies the coexistence of stable equilibria with stable and unstable periodic solutions in a narrow range of volumetric flow rates for Q in [0.5215, 0.5364]. This range is more clearly seen in Fig. 9, where the period of both stable (solid) and unstable (dashed) limit cycles is plotted versus the volumetric flow rate. Transient simulations of slightly perturbed unstable limit cycles can lead to either stable steady states or stable periodic solutions as shown in Fig. 10, where the two possibilities when starting with the unstable limit cycle for Q=0.530 are illustrated: when the volumetric flow rate is slightly perturbed to Q=0.531, the solution evolves to the corresponding steady state (Fig. 10a); when Q is set to 0.529 instead, the system reaches the corresponding stable periodic solution (Fig. 10b). The trajectories of the above solutions are shown in Fig. 11, where the phase portraits of the slip velocity versus the viscoelastic stress $T_1^{xy}\big|_{y=0.5}$ are plotted.

## 5. Conclusions

The timestepper approach has been employed for the stability and bifurcation analysis of the plane Poiseuille flow of an Oldroyd-B fluid with nonmonotonic slip at the wall. The effectiveness of the method in determining the regimes of instability has been demonstrated for given values of the elasticity number and the Newtonian viscosity. Unstable steady-states exist in the negative-slope regime of the flow curve in a range of volumetric flow rates defined by two Hopf bifurcation points. The existence of both stable and unstable periodic solution branches has been determined by means of a pseudo arc-length continuation algorithm treating the imposed volumetric flow rate as the bifurcation parameter.
In the linearly unstable regime between the two Hopf points, there exist only stable periodic solutions, i..e. self-sustained oscillations of the pressure gradient, the velocity, and the shear stress at fixed volumetric flow rate. The most interesting finding of this work is the coexistence of stable equilibria and both stable and unstable periodic solutions in a narrow range of flow rates beyond the linearly unstable regime, i.e. the right Hopf bifurcation point is a subcritical one.

**Acknowledgements**
The work of C.I.S. was supported by the European Social Fund and National Resources - (EPEAEK II) –Pythagoras II. M.E.K. and A.G.B. acknowledge partial support by the State Scholarships Foundation of Greece and the National Technical University of Athens through the Basic Research Program ''Protagoras''.

# Figure Captions

**Figure 1** Computation of fixed points using the timestepper approach and GMRES.

**Figure 2** Geometry and boundary conditions in plane Poiseuille flow of an Oldroyd-B fluid with slip at the wall.

**Figure 3** Flow curve for $A_1=1$, $A_2=15$ and $A_3=100$.

**Figure 4** Critical eigenvalues of the steady state Jacobian matrix at **(a)** Q=0.413 and **(b)** Q=0.414. The transition of the two complex eigenvalues from the left real plane to the right plane indicates the existence of a supercritical Hopf point.

**Figure 5** Critical eigenvalues of the steady state Jacobian matrix at **(a)** Q=0.521 and **(b)** Q=0.522. The transition of two complex eigenvalues from the right real plane to the left one indicates the existence of a subcritical Hopf point.

**Figure 6** Evolution of **(a)** the slip velocity and **(b)** the pressure gradient when the unstable steady-state solution at Q=0.449 is perturbed by setting Q=0.45.

**Figure 7** Bifurcation diagram of stable / unstable steady and periodic solutions for the plane Poiseuille flow of an Oldroyd-B fluid with non monotonic slip at the wall. The solid line corresponds to stable and the dashed line to unstable steady states. The solid line with black circles corresponds to stable limit cycles, while the dashed line with open circles corresponds to unstable ones. In the case of limit cycles, the maximum value of the pressure gradient during one period is plotted. The inset is a blow up of the bifurcation diagram in the neighborhood of the turning point of limit cycles.

**Figure 8** Dominant Floquet multipliers for the determination of the periodic solution stability: **(a)** the leading Floquet multiplier has modulus less than one indicating a stable periodic solution; **(b)** the leading multiplier has modulus greater than one which means that the periodic solution is unstable. There is always a trivial Floquet multiplier equal to one.

**Figure 9** The period *T* of limit cycles as a function of the volumetric flow rate. The solid line corresponds to stable and the dashed one to unstable limit cycles.

**Figure 10** Evolution of the pressure gradient and the slip velocity starting from the unstable periodic solution for *Q*=0.530 (dashed line) and perturbing the volumetric flow rate to **(a)** *Q*=0.531 (the system evolves towards a stable steady state) and **(b)** Q=0.529 (the system evolves towards a stable limit cycle)

**Figure 11** Phase portraits of slip velocity versus the viscoelastic stress component $T_1^{xy}\big|_{y=0.5}$. The application of a small perturbation to the unstable limit cycle (dashed line) obtained at *Q*=0.530 leads to: **(a)** a stable steady state when the volumetric flow rate is perturbed to 0.531 and **(b)** a stable periodic solution when the volumetric flow rate is set to 0.529.

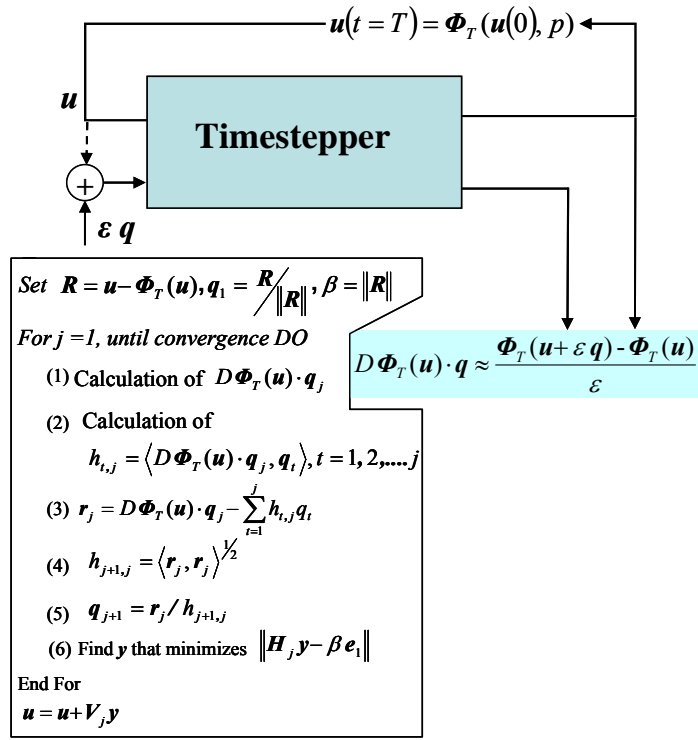

**Figure 1**

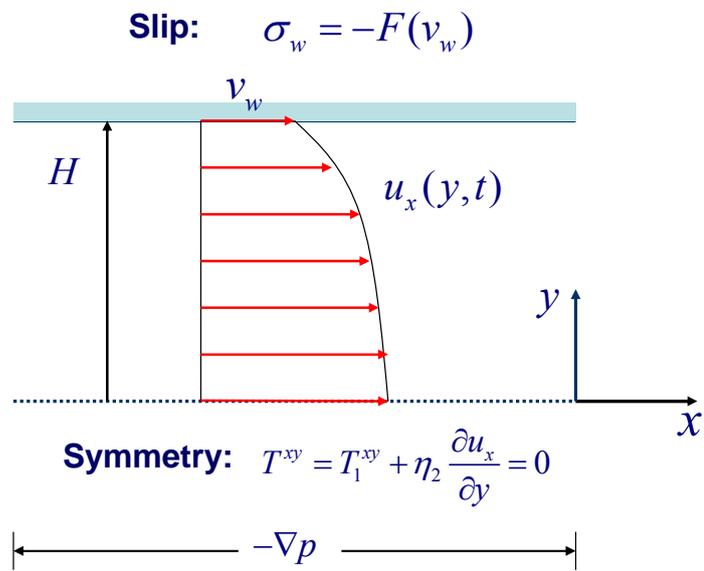

**Figure 2**

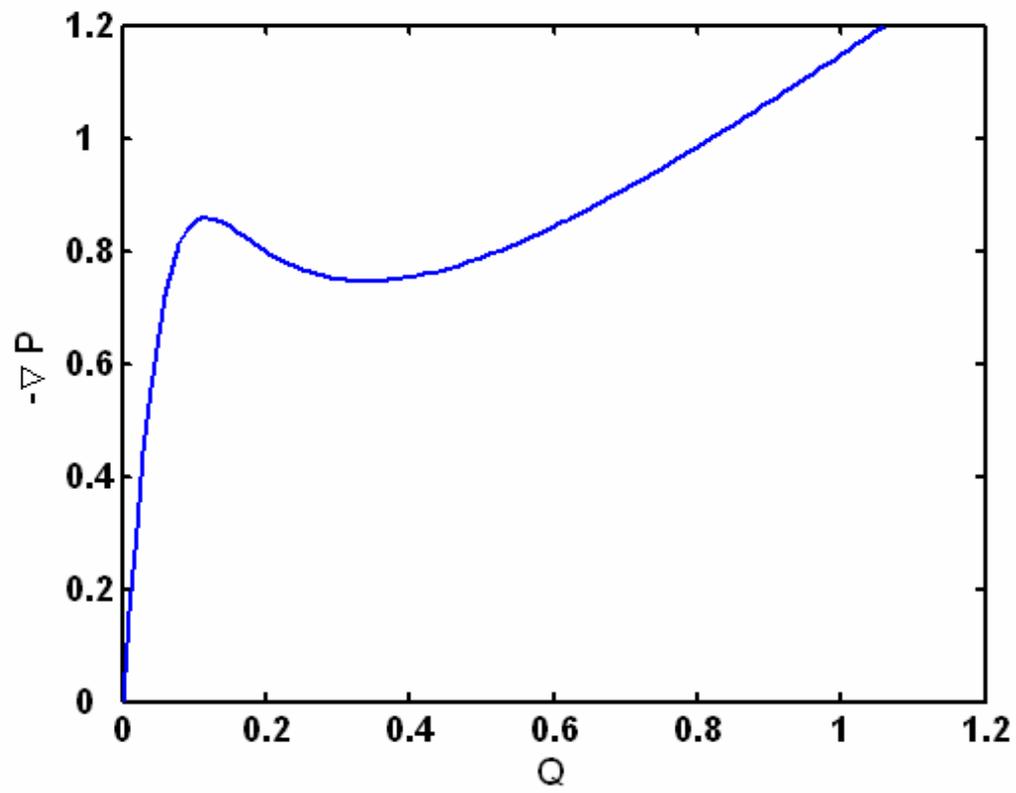

**Figure 3**

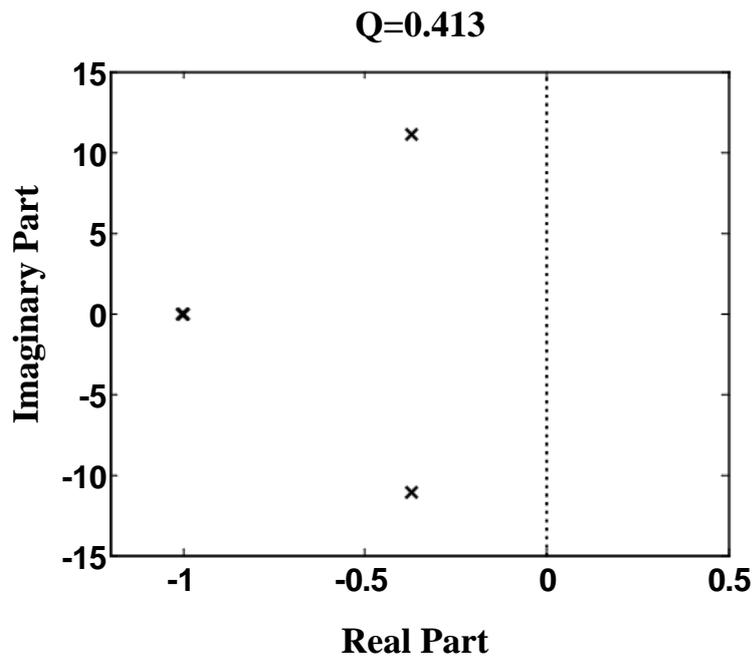

(a)

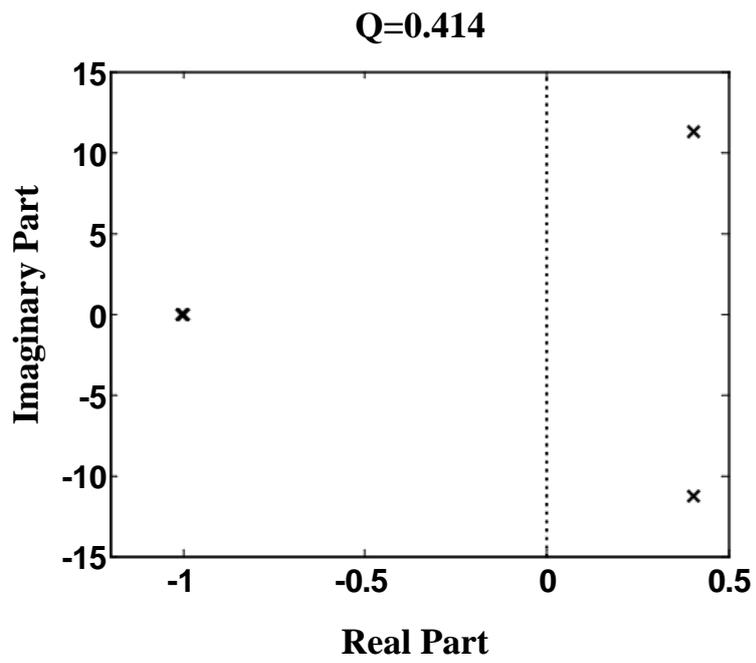

(b)

**Figure 4**

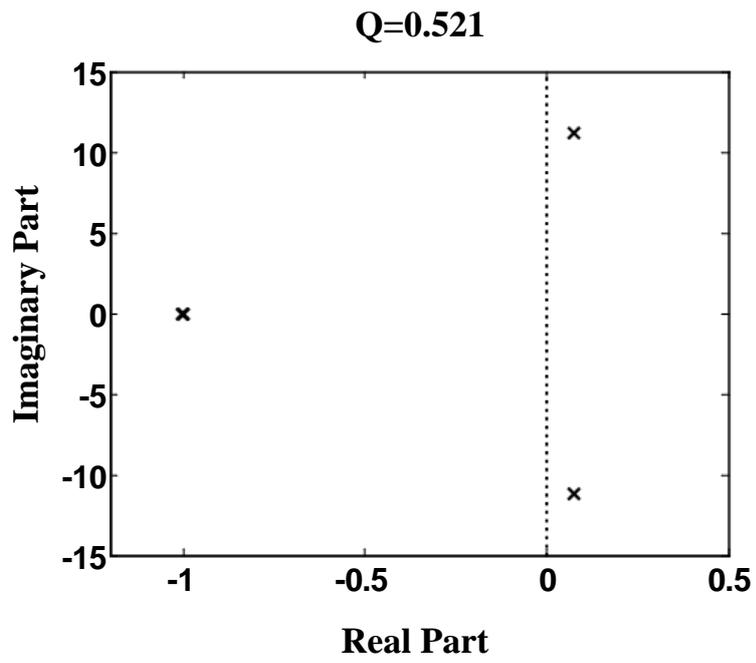

(a)

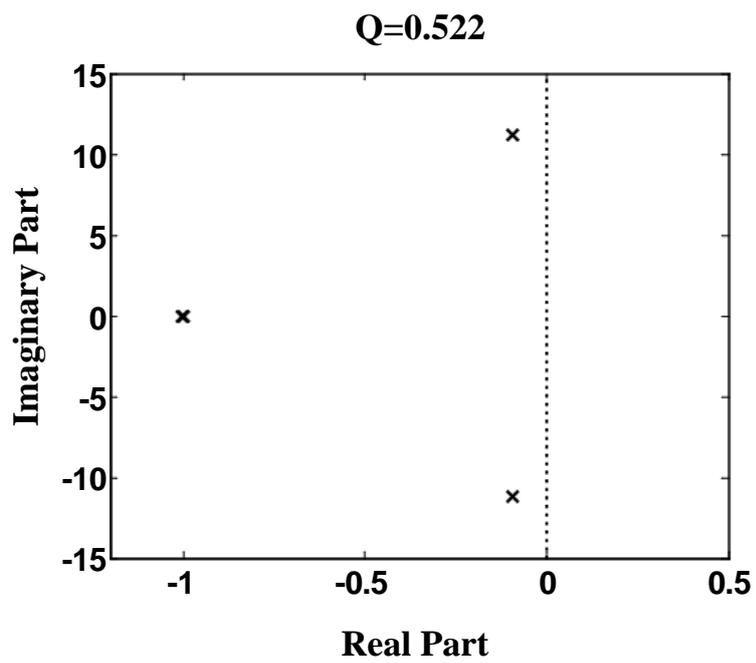

(b)

**Figure 5**

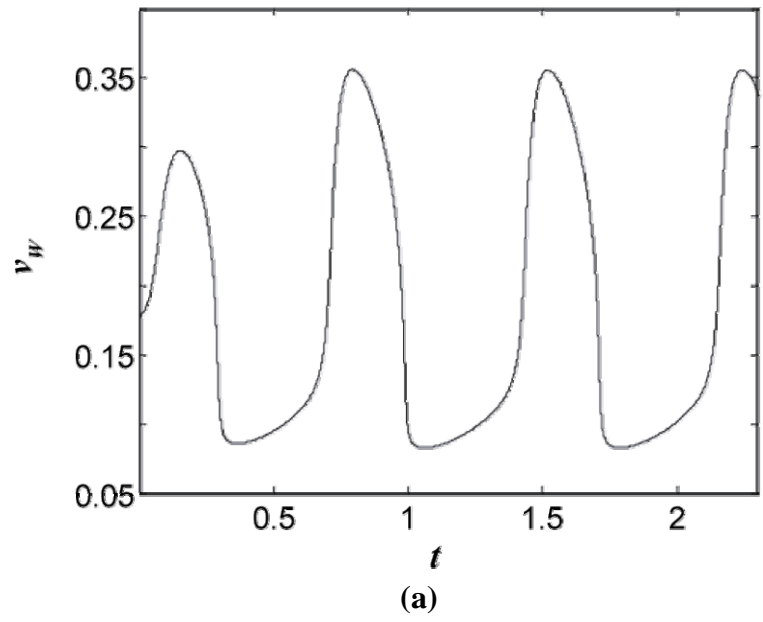

(a)

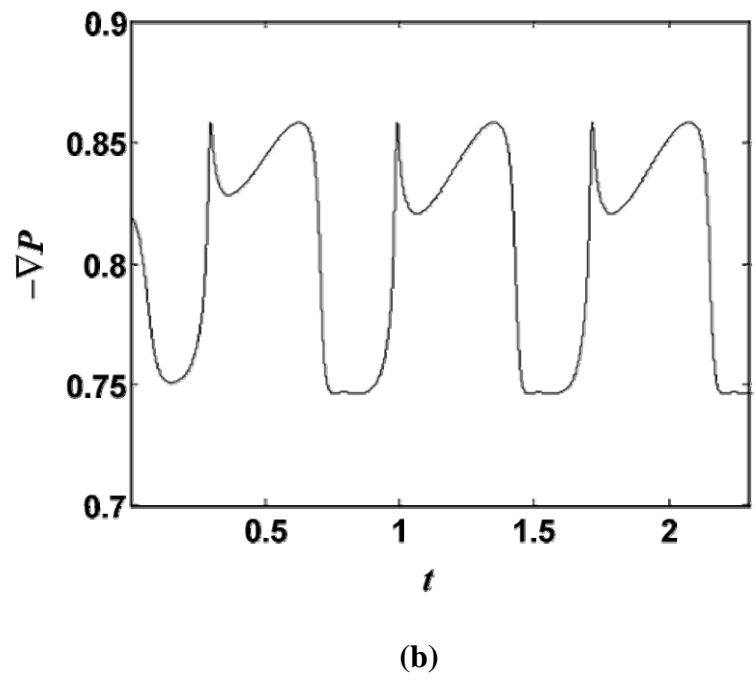

(b)

**Figure 6**

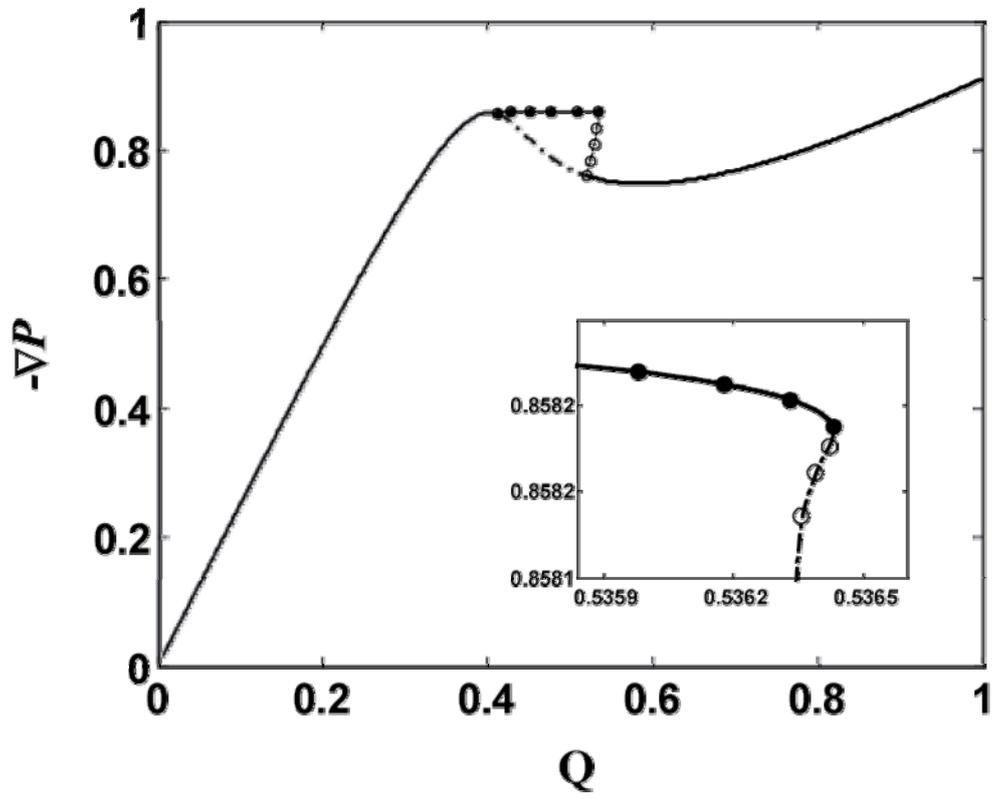

**Figure 7**

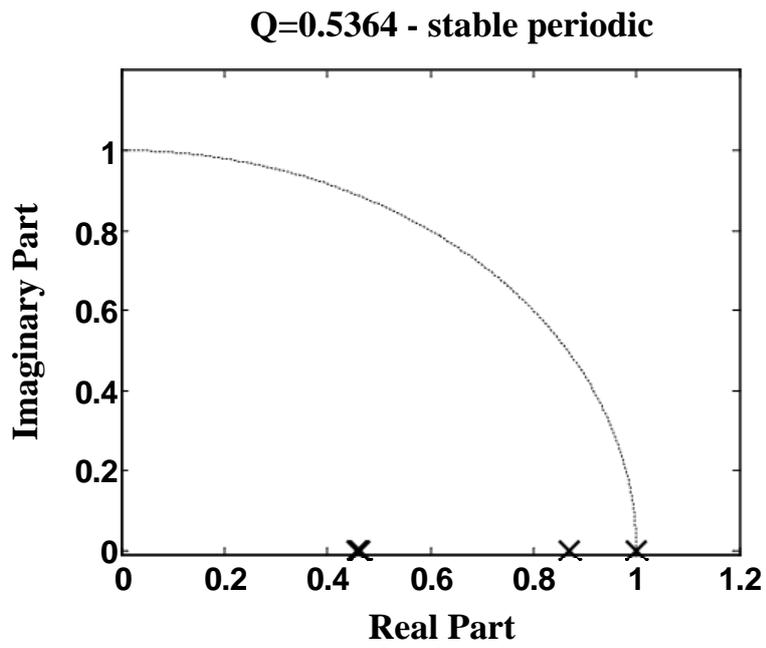

(a)

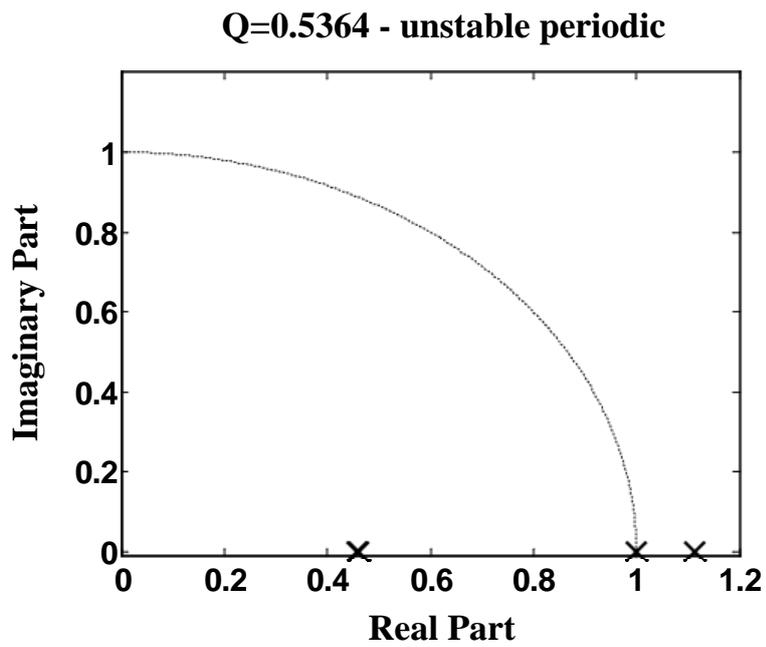

(b)

**Figure 8**

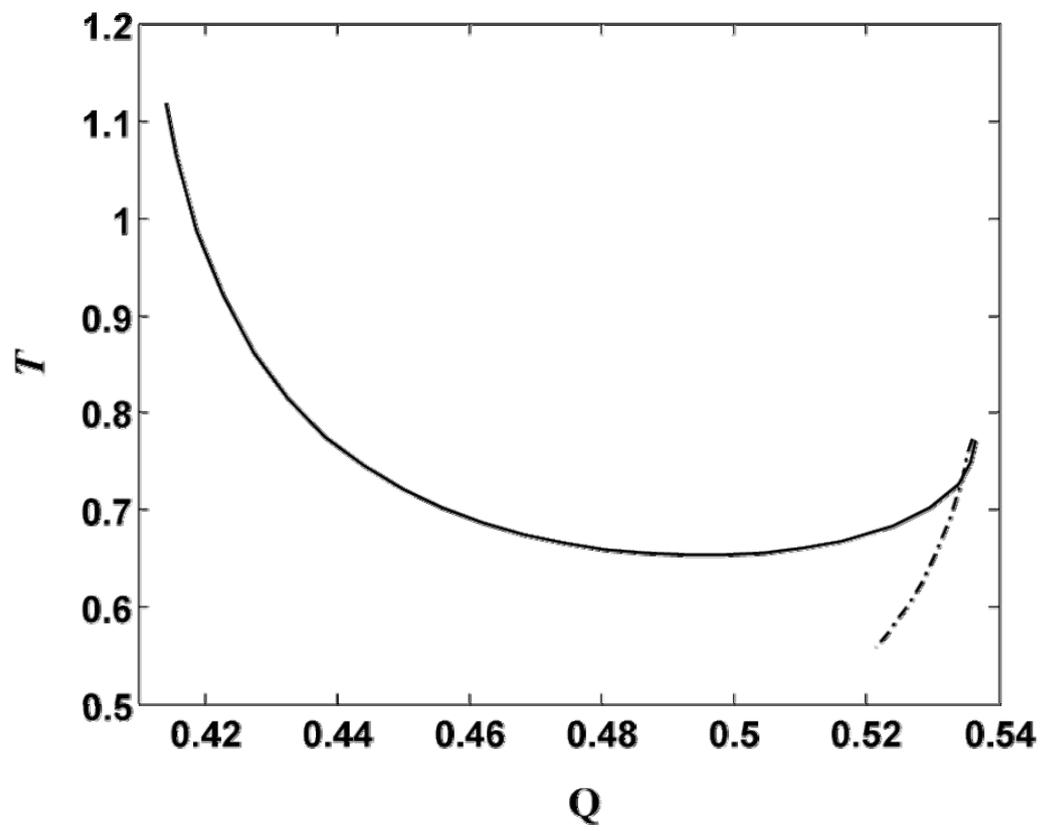

**Figure 9**

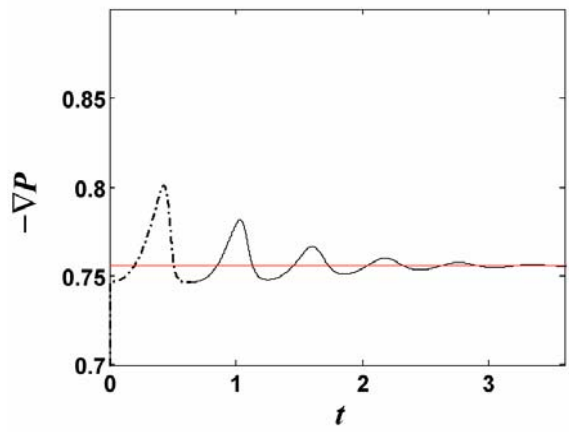 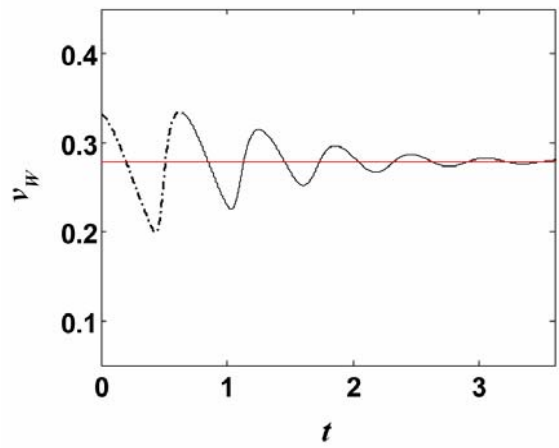

**(a)**

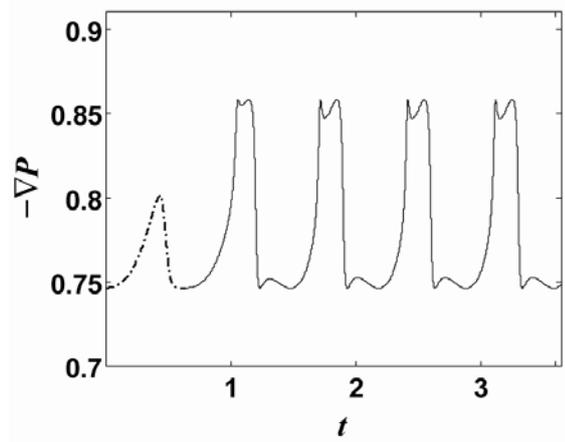 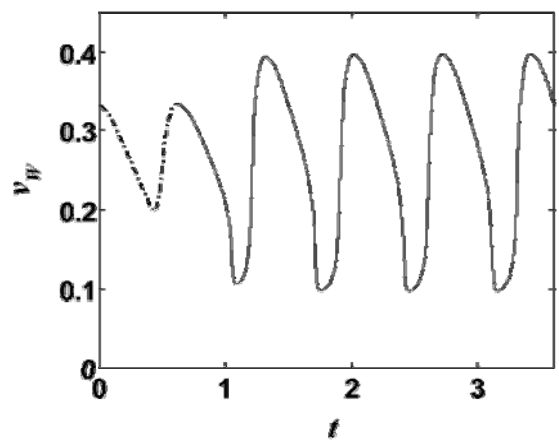

**(b)**

**Figure 10**

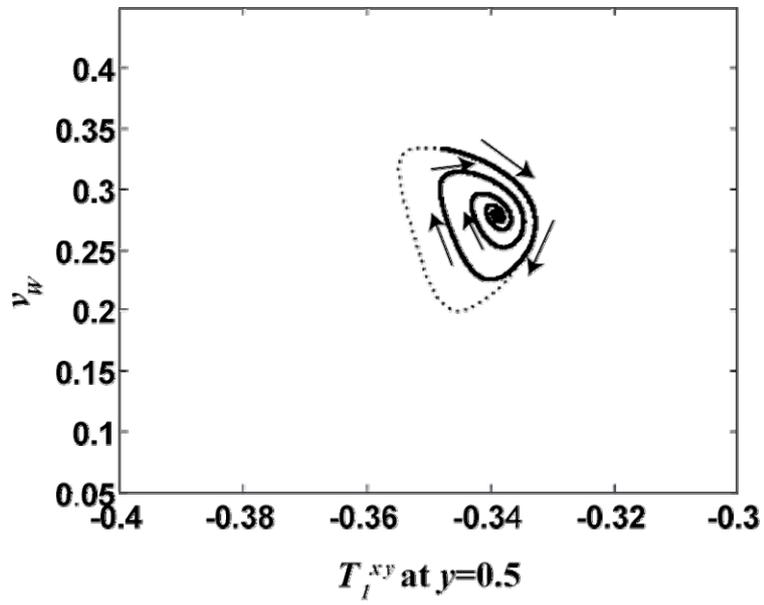

(a)

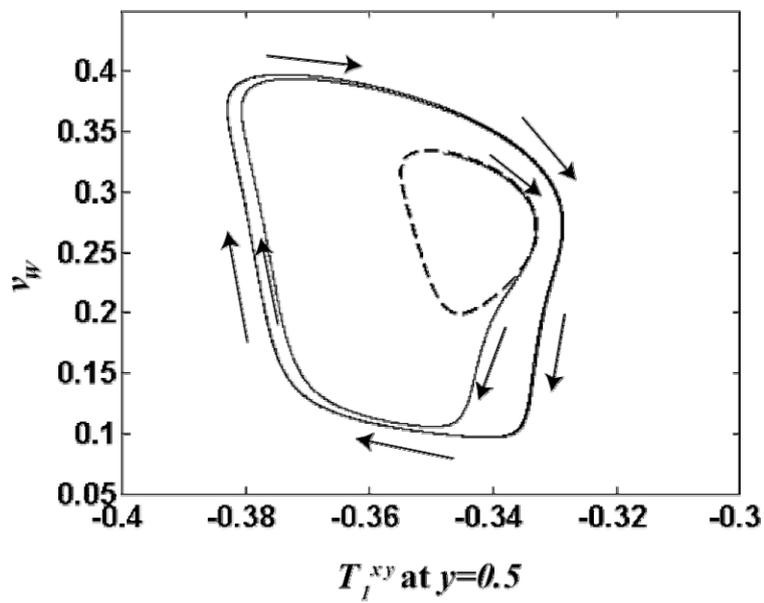

(b)

**Figure 11**